\documentclass[11pt,a4paper,reqno]{amsart}

\usepackage[utf8]{inputenc}
\usepackage[english]{babel}
\usepackage{amssymb,mathrsfs}
\usepackage{hyperref}

\calclayout
\hypersetup{colorlinks=true,linkcolor=blue,citecolor=blue,urlcolor=blue}

\theoremstyle{definition}
\newtheorem{theorem}{Theorem}[section]
\newtheorem{lemma}[theorem]{Lemma}
\newtheorem{proposition}[theorem]{Proposition}
\newtheorem{remark}[theorem]{Remark}
\numberwithin{equation}{section}

\newcommand{\cP}{\mathcal{P}}
\newcommand{\R}{\mathbb{R}}
\newcommand{\CC}{\mathbb{C}}
\newcommand{\cH}{\mathcal{H}}
\newcommand{\E}{\mathcal{E}}
\newcommand{\F}{\mathcal{F}}
\newcommand{\cS}{\mathcal{S}}
\newcommand{\Bin}{\mathcal{B}^{\text{in}}}
\newcommand{\Bout}{\mathcal{B}^{\text{out}}}
\newcommand{\We}{\mathrm{We}}
\newcommand{\cI}{\mathcal{I}}
\newcommand{\dd}{\mathop{}\!\mathrm{d}}
\newcommand{\ii}{\mathrm{i}}
\newcommand{\ee}{\mathrm{e}}
\DeclareMathOperator{\capa}{cap}
\newcommand{\abs}[1]{\left| #1 \right|}
\def\N{\mathbb{N}}
\def\Z{\mathbb{Z}}
\def\norm#1{\left\| #1 \right\|}
\newcommand{\C}{\operatorname{C}}  
\newcommand{\Hardy}{\mathcal{H}}  



\begin{document}
\allowdisplaybreaks

\title{Global rigidity of two-dimensional bubbles}

\date{\today}

\author{Lukas Niebel}
\email{lukas.niebel@uni-muenster.de }
\address[Lukas Niebel]{Institut f\"ur Analysis und Numerik, Universit\"at M\"unster\\
	Orl\'eans-Ring 10, 48149 M\"unster, Germany.}

\begin{abstract}
	We study bubbles, which are stationary hollow vortices with surface tension, in two dimensions. Such objects solve an overdetermined elliptic free boundary problem in an exterior domain, with an additional boundary condition involving curvature and the Neumann trace. Motivated by a conjecture of Crowdy and Wegmann, we prove that every solution prescribed by a Jordan curve is circular for $0 \le \mathrm{We} \le 3$, where $3$ is the first bifurcation value. This range is sharp: the bifurcation branch at $\mathrm{We}=3$ contains non-circular solutions with $\mathrm{We} > 3$ arbitrarily close to $3$. This elliptic problem describes critical points of the sum of the perimeter and the logarithmic potential energy of compact sets. We prove an area-constrained perimeter--capacity inequality and show that, in the associated convexity-constrained variational problem, the unit disc is the unique minimiser up to translation if and only if $0 \le \mathrm{We} \le 3$.
\end{abstract}

\maketitle

\section{Introduction}
We study the \emph{free boundary Euler equations} with \emph{surface tension} in two dimensions:
\begin{equation*}
	\begin{aligned}
		\rho(\partial_t U + (U \cdot \nabla) U) + \nabla P & = 0             &  & \quad\mbox{ for } t \in \R, \, x \in \Bout(t) \\
		\nabla \cdot U                                     & = 0             &  & \quad\mbox{ for } t \in \R, \, x \in \Bout(t) \\
		-P                                                 & = \sigma \kappa &  & \quad\mbox{ for } t \in \R, \, x \in \cS(t)   \\
		U \cdot n                                          & = V_n           &  & \quad\mbox{ for } t \in \R, \, x \in \cS(t).  \\
	\end{aligned}
\end{equation*}
Here $U$ denotes the fluid velocity, $P$ the hydrodynamic pressure, and $\rho>0$ the constant mass density of the fluid phase.
The first two relations express momentum conservation away from the interface $\cS(t)$ and incompressibility.
Because the free interface is a moving surface carried by the fluid, its normal speed $V_n$ (i.e. $\partial_t \cS \cdot n = V_n$) must equal the fluid's normal velocity $U\cdot n$ at the interface.
Motivated by bubble dynamics, we assume that the inner phase occupies a bounded and connected domain $\Bin(t)$ and is dynamically inert, with spatially constant pressure $P_{\mathrm{in}}(t)$.
The surrounding fluid is $\Bout(t) = \R^2\setminus \overline{\Bin(t)}$, so that $\cS(t)=\partial \Bin(t)=\partial \Bout(t)$.
Taking $P_{\mathrm{in}}(t)$ as the pressure reference, the third relation is the \emph{Young--Laplace law}, where $\sigma>0$ is the surface tension and $\kappa$ the signed curvature.
We call these objects bubbles, with a submerged air column as the model case.

We choose the signed curvature $\kappa$ to be positive for convex $\Bin$, and normalise it so that $\kappa=1$ when $\Bin$ is the unit disc.
The unit normal vector $n$ points from the inner phase into the outer one.

\subsection*{Stationary bubbles}
To find a \emph{stationary bubble}, also called a \emph{stationary hollow vortex solution} $U = U(x)$ to the free boundary Euler equations with surface tension and irrotational outer flow, we have to solve the following \emph{overdetermined elliptic free boundary value problem} for the stream function $\psi = \psi(x)$ with $U = \nabla^\perp \psi = (-\partial_2 \psi , \partial_1 \psi)$. The logarithmic asymptotics imposed below ensure that the fluid velocity vanishes at infinity. Since the stream function is defined up to an additive constant, we normalise it by requiring $\psi=0$ on the free boundary.

Let $R>0$, $\rho>0$, $\alpha \in \R$ and $\sigma \ge 0$.
The task is to find a $\C^{1,1}$ closed Jordan curve $\cS$ partitioning $\R^2$ into an interior domain
$\Bin$ with $\abs{\Bin} = \pi R^2$ and an exterior domain $\Bout$, and a solution $\psi \colon \Bout \to \R$ to
\begin{equation} \label{eq:hollowvortex}
	\begin{aligned}
		-\Delta \psi                                     & = 0                          &  & \quad \mbox{ in } \Bout             \\
		\psi                                             & = 0                          &  & \quad \mbox{ on } \cS               \\
		\psi(x)                                          & = \alpha \log \abs{x} + O(1) &  & \quad \mbox{ as } \abs{x}\to \infty \\
		-\frac{\rho}{2}(\partial_n \psi)^2+\sigma \kappa & = \lambda                    &  & \quad \mbox{ on } \cS,
	\end{aligned}
\end{equation}
for some $\lambda \in \R$.
The last equation is called the jump equation.
By integrating over large circles anticlockwise, we observe that $2\pi\alpha$ coincides with the total circulation.
To rewrite the jump condition for the pressure law, we used Bernoulli's law for steady flows (see \cite[Chapter 1.9]{MR1217252}).

We focus on the capillary regime $\sigma>0$ with non-zero circulation $\alpha\neq0$. The key quantity to describe the solutions of the overdetermined free boundary value problem is the \emph{Weber number}
\begin{equation*}
	\We = \frac{\rho\alpha^2}{\sigma R}.
\end{equation*}

\subsection*{The overdetermined elliptic free boundary value problem}
We nondimensionalise the problem as follows: find $\We \ge 0$, $\lambda \in \R$, a $\C^{1,1}$ closed Jordan curve $\cS$ partitioning $\R^2$ into an interior domain $\Bin$ with $\abs{\Bin} = \pi $ and an exterior domain $\Bout$, and a solution $\psi \colon \Bout \to \R$ to
\begin{equation} \label{eq:ofbvp}
	\begin{aligned}
		-\Delta \psi                                & = 0                   &  & \quad \mbox{ in } \Bout             \\
		\psi                                        & = 0                   &  & \quad \mbox{ on } \cS               \\
		\psi(x)                                     & = \log \abs{x} + O(1) &  & \quad \mbox{ as } \abs{x}\to \infty \\
		-\frac{1}{2}\We (\partial_n \psi)^2+ \kappa & = \lambda             &  & \quad \mbox{ on } \cS.
	\end{aligned}
\end{equation}

Let us first explain the known solutions to \eqref{eq:ofbvp}.
Every unit circle centred at $x_0 \in \R^2$ is a solution to \eqref{eq:ofbvp}.
In fact, the function $\psi = \log \abs{x-x_0}$ solves the free boundary value problem and satisfies the jump condition
\begin{equation*}
	-\frac{1}{2} \We + 1 = \lambda.
\end{equation*}
The first non-circular solutions were constructed by Crowdy \cite{Crowdy1999}, who found an explicit family of two-fold symmetric bubbles using rational conformal maps. Wegmann and Crowdy \cite{Wegmann2000-ql} subsequently generalised this construction and obtained a branch of $m$-fold symmetric bubbles for every integer $m\ge 2$. The $m$-fold branch bifurcates from the circular solution at $\We=m+1$.
Crowdy and Wegmann conjecture ``that the bubble is circular for most values of circulation [equivalently of $\We$] with the exception of an infinite series of discrete values'' in \cite{Wegmann2000-ql}.
Here the discrete values refer to the points at which
non-circular branches bifurcate from the circular family.
The main goal of this article is to prove \emph{global rigidity} of the circular solution for $0 \le \We \le 3$, which further supports this conjecture as $\We =3$ is the smallest bifurcation value. We obtain the following

\begin{theorem}\label{thm:globalWe3}
	Let $(\cS,\psi,\We,\lambda)$ be a solution to \eqref{eq:ofbvp}. If $\cS$ is not a circle, then $\We>3$.
\end{theorem}

No geometric assumption such as convexity, symmetry or being close to a circle is imposed.

We prove Theorem \ref{thm:globalWe3} in Section \ref{sec:globalWe3}.
We first use an exterior conformal map to rewrite the free boundary problem as a problem on the unit circle.
Applying Fourier estimates to an analytic square root of the derivative of the conformal map and using the area constraint, we obtain a lower bound for $\We$ in terms of the perimeter.
The isoperimetric inequality then implies that $\We>3$ unless $\cS$ is a circle.

We also prove the weaker conclusion $\We>2$ in Section \ref{sec:global} by
combining geometric identities with identities for solutions of the elliptic
PDE, such as the Pohozaev identity \cite{MR192184}. We believe that these ideas
may be useful in higher dimensions.

The result is sharp because the two-fold Crowdy--Wegmann branch contains non-circular solutions with $\We>3$ arbitrarily close to $3$ \cite{Wegmann2000-ql}.
By closely investigating the Crowdy--Wegmann branches we prove in Section \ref{sec:cw} that both a circular
and a non-circular solution exist for every $\We \ge 10$. Numerical evidence for the remaining range
$3<\We<10$ is presented in Remark \ref{rem:noncircular}.

\subsection*{Variational bubbles}

A related \emph{variational problem} will be investigated in Section \ref{sec:var}, where we obtain an analogous global rigidity result for its minimisers.
Equation \eqref{eq:ofbvp} has a variational structure analogous to non-local
charged-drop models \cite{Fontelos2004-ym,Goldman2015-aj,MR3785601}, with the
logarithmic potential playing the role of the electrostatic potential.

For compact sets $E \subset \R^2$ of positive area and finite perimeter, we consider the functional
\begin{equation*}
	\F_{\We}(E) = \We \pi \ \cI(E) + \cP(E)
\end{equation*}
with the logarithmic potential energy
\begin{equation*}
	\cI(E) = \inf_{\substack{\mu \text{ probability measure} \\\mu(E) = 1}} \left\{- \int_{\R^2} \int_{\R^2}  \log \abs{x-y} \dd \mu(y) \dd \mu(x) \right\},
\end{equation*}
and where $\cP(E)$ denotes the perimeter. For smooth critical points,
the correspondence with the overdetermined elliptic free boundary value problem
\eqref{eq:ofbvp} follows from Hadamard's variational formula for the exterior
Green function with pole at infinity; see
\cite{MR5020640}.
We refer to \cite{MR350027} for background on logarithmic potential theory.
We introduce the logarithmic capacity as $\capa(E)=\ee^{-\cI(E)}$.

For $\We>0$, the variational problem
\begin{equation*}
	\inf_{\substack{E \subset \R^2\\ E \text{ compact}\\ \abs{E} = \pi}} \F_{\We}(E)
\end{equation*}
is ill-posed, as can be seen by evaluating it on two discs and sending the distance between their centres to $\infty$.
To overcome this, one could impose a further geometric constraint on $E$, such as convexity or a $\delta$-ball condition; see \cite{Goldman2015-aj,MR3785601}.
We restrict the minimisation problem to convex sets, since, for general Weber numbers, it is not known whether the infimum over closures of bounded domains with boundaries being rectifiable Jordan curves is attained.
For the variational problem of minimising this functional over all compact and convex sets
\begin{equation} \label{eq:infconvF}
	\inf_{\substack{E \subset \R^2\\E \text{ compact and convex} \\ \abs{E} = \pi}} \F_{\We}(E)
\end{equation}
global rigidity of the unit disc for small Weber numbers $\We \le \We_0$, for some small, quantitative but implicit constant $\We_0$, and failure of minimality for sufficiently large $\We$ were proved in \cite{MR3785601}. The small-Weber-number argument depends, among other ingredients, on the sharp quantitative isoperimetric inequality due to Fusco, Maggi and Pratelli \cite{MR2456887}.
We improve this classification to the sharp range $0 \le \We \le 3$ by proving a new inequality relating the perimeter and the logarithmic capacity of sets bounded by a Jordan curve of fixed area. Our second main result is the following

\begin{theorem}
	\label{thm:exact-threshold}
	In the class of convex and compact sets $E\subset\R^2$ with $\abs{E}=\pi$,
	the unit disc is the unique minimiser of $\F_{\We}$, up to translation,
	if and only if $0\le \We\le3$.
\end{theorem}

The key ingredient for the proof of this theorem is a new perimeter--capacity inequality
\[
	\cP(E)-2\pi
	\ge
	2\pi\bigl(\sqrt{5\capa(E)^2-1}-\capa(E)-1\bigr),
\]
for compact sets $E\subset\R^2$ whose boundary is a rectifiable Jordan curve and which satisfy $\abs{E}=\pi$.
The area constraint makes
this a refinement of the classical perimeter--capacity inequality
$2\pi\capa(E)\le\cP(E)$, which goes back to P\'olya--Schiffer and
P\'olya--Szeg\H{o} \cite{MR105507,MR396134,PolyaSzego1951}.

The proof uses the normalised exterior conformal map of $E$. Gr\"onwall's
area formula and a boundary representation of the perimeter express the
area constraint and the perimeter in terms of a Bergman norm and a Hardy
norm, respectively. Then, we take an analytic square root of the normalised
derivative and apply Burbea's sharp product inequality
\cite{MR882113} to deduce the stated estimate.

Formally, the Euler--Lagrange equation is the overdetermined free boundary value problem \eqref{eq:ofbvp}. For convex sets, however, the admissible perturbations are restricted, and the last equation in \eqref{eq:ofbvp} follows only on the strictly convex part of the boundary; see \cite[Proposition 1]{MR2927624}. It is therefore unclear whether minimisers of \eqref{eq:infconvF} solve the full overdetermined problem, so Theorem \ref{thm:globalWe3} does not imply Theorem \ref{thm:exact-threshold}.

\subsection*{Literature}
The two limiting cases of the problem are classical rigidity problems.
For $\We = 0$, the jump equation prescribes constant curvature, and hence $\cS$ must be a unit circle up to translation.
In the formal limit $\We = \infty$, surface tension vanishes and \eqref{eq:hollowvortex} reduces to the exterior overdetermined problem with $\sigma = 0$, $\alpha = 1$ and $R = 1$.
Serrin proved the corresponding rigidity theorem for bounded domains \cite{MR333220}, while Reichel proved rigidity of the unit circle in the exterior setting \cite[Theorem 2]{MR1463801}.
Our argument is inspired by Weinberger's proof of Serrin's theorem \cite{MR333221}.

The exact families recalled above originate in Crowdy's rational conformal maps for two-fold symmetric bubbles \cite{Crowdy1999}.
Wegmann and Crowdy subsequently obtained $m$-fold symmetric families of solutions to \eqref{eq:ofbvp} for every integer $m\ge 2$.
They numerically
identified their first pinching thresholds for $m=2,\ldots,7$
\cite[p.~2136 and Sec.~5]{Wegmann2000-ql}.
Hayford and Wang \cite{HayfordWang2022} later revisited these families in an
equivalent electrified-droplet formulation and obtained the exact
pinching thresholds for the cases $m=2,3,4$.
The conformal mapping framework was developed further in \cite{MR1856629}.
The bifurcation values $\We = m+1$ form a discrete sequence, whereas each nontrivial branch contains a continuum of Weber numbers.
The linear analysis also indicates local rigidity of the unit circle whenever $\We\notin\{n\in\N \; ; \; n\ge3\}$ as remarked in \cite{Wegmann2000-ql}. This can be made rigorous by a standard implicit-function argument, for instance using the framework of \cite{mns_bubbles_25,Fontelos2004-ym}; since this result is not used below, we omit the details.

Crowdy, Nelson and Krishnamurthy \cite{MR4226915} constructed rotating hollow vortices without surface tension and called them $H$-states.
More generally, Crowdy and Tanveer \cite{MR4965079} showed in several two-dimensional settings that constant-vorticity flows without surface tension and irrotational flows with surface tension reduce to the same mathematical free boundary problem.
Thus, the stationary capillary bubbles considered here and the rotating $H$-states share the same free boundary profiles, although their physical equations are different.
We refer to \cite{MR3224254} for the relation between hollow vortices, capillary water waves and double quadrature domains.

The situation in the case of zero surface tension is quite different when the geometry or topology is broadened.
Exceptional domains and related overdetermined problems were studied in \cite{MR2794602,MR3095114}.
Traizet classified planar solutions and related hollow vortices to minimal surfaces in \cite{MR3192039,MR3455335}.
Quasi-exceptional domains were studied by Eremenko and Lundberg \cite{MR3366032}.
Moreover, non-degenerate point-vortex equilibria were desingularised into translating, rotating or stationary hollow vortices in \cite{MR5037134}.

The variational problem is related to non-local models for charged liquid drops \cite{Fontelos2004-ym,Goldman2015-aj,MR3785601}.
Crowdy \cite{Crowdy2015ChargedDewetting} identified a further connection between the rational two-fold profile from \cite{Crowdy1999} and static charged conducting droplets on a substrate.

We finally discuss a related three-dimensional result \cite{mns_bubbles_25}.
Close-to-spherical steady bubbles and drops with a vortex core were constructed in the three-dimensional setting.
The three-dimensional axisymmetric problem with irrotational flow around the bubble is fundamentally different from its two-dimensional counterpart.
For $\We = 0$, the sphere is the unique solution.
For small $\We > 0$, the constructed solutions form the unique branch emerging from the sphere and contain all close-to-spherical solutions in the class considered there.
To our knowledge, no further axisymmetric irrotational solutions with the topology of a sphere are known in this setting.
It remains open whether every axisymmetric irrotational solution with sufficiently small Weber number belongs to this branch.

\subsection*{Acknowledgements}
The author thanks Bj\"orn Gebhard, Yuanjiang Han and Christian Seis for fruitful discussions. Moreover, the author is very grateful to the referees who helped to significantly improve the presentation of the manuscript.

Lukas Niebel is funded by SNSF Starting Grant
TMSGI2\textunderscore226018 and by the Deutsche Forschungsgemeinschaft (DFG,
German Research Foundation) under Germany's Excellence Strategy
EXC 2044/2--390685587, Mathematics M\"unster:
Dynamics--Geometry--Structure and grant no.~531098047.

\vfill

\subsection*{Data availability statement} Not applicable.

\subsection*{Conflict of interest} The author declares that he has no conflict of interest.

\section{Pohozaev and geometric identities: a preliminary rigidity result}
\label{sec:global}

This section proves the following weaker version of Theorem \ref{thm:globalWe3}.
\begin{proposition}\label{prop:globalWe2}
	Let $(\cS,\psi,\We,\lambda)$ be a solution to \eqref{eq:ofbvp}. If $\cS$ is not a circle, then $\We>2$.
\end{proposition}

Fix such a solution.
First, we note that $\cS$ has to be smooth, see \cite{KochLeoni}.
We recall that $\kappa$ denotes the signed curvature with respect to the unit normal
$n$ pointing from $\Bin$ into
$\Bout$. Let $L:=\cH^1(\cS)=\cP(\Bin)>0$ and
$\gamma\colon\R/(L\Z)\to \cS$ be an arclength
parametrisation oriented anticlockwise. Set
$n(s):=n(\gamma(s))$, $\kappa(s):=\kappa(\gamma(s))$, and
$\tau(s):=\gamma'(s)$. Then
\[
	\tau'(s)=\gamma''(s)=-\kappa(s)n(s).
\]
Writing $x=r(\cos\theta,\sin\theta)$, for sufficiently large $r$ the exterior Fourier
expansion of $\psi$, with $A_0,A_k,B_k\in\R$ for $k\ge1$, is
\[
	\psi(r,\theta)
	=\log r+A_0+
	\sum_{k=1}^{\infty}r^{-k}
	\bigl(A_k\cos(k\theta)+B_k\sin(k\theta)\bigr).
\]
Termwise differentiation gives
\begin{equation}
	\partial_r\psi
	=\frac1r+O\left(\frac1{r^2}\right),
	\qquad
	\abs{\nabla\psi}^2
	=\frac1{r^2}+O\left(\frac1{r^3}\right)
	\quad\text{as }r\to\infty,\ \text{uniformly in }\theta.
	\label{eq:far-field-derivatives}
\end{equation}

We collect two simple identities.
The divergence theorem, harmonicity of $\psi$, and
\eqref{eq:far-field-derivatives} give
\begin{equation} \label{eq:fluxL1}
	\int_{\cS} \partial_n \psi \dd \cH^1 = 2\pi.
\end{equation}
Since $\psi=0$ on $\cS$, its tangential derivative vanishes there, and hence
$\abs{\nabla\psi}=\abs{\partial_n\psi}$ on $\cS$.

Integrating the jump condition and using the Gau\ss--Bonnet theorem yields
\begin{equation} \label{eq:jumpGB}
	- \frac{\We}{2} \int_\cS (\partial_n \psi)^2 \dd \cH^1 + 2\pi  = \lambda \cP(\Bin).
\end{equation}

Next, we use the Pohozaev identity \cite{MR192184} to derive a key identity. For the convenience of the reader we recall that if $u\colon \Omega \subset \R^d \to \R$, $d \in \N$, is harmonic $\Delta u = 0$, then the vector field
\begin{equation} \label{eq:pohozaev-vector}
	X(x):=\bigl(x\cdot\nabla u(x)\bigr)\nabla u(x)
	-\frac12\lvert\nabla u(x)\rvert^2x,
\end{equation}
satisfies
\begin{equation} \label{eq:pohozaev}
	\nabla\cdot X(x)
	=\frac{2-d}{2}\lvert\nabla u(x)\rvert^2.
\end{equation}

\begin{lemma} \label{lem:pohozaev}
	We have
	\begin{equation*}
		\int_{\cS}(x \cdot n) (\partial_n \psi)^2 \dd \cH^1 = 2\pi .
	\end{equation*}
\end{lemma}

\begin{proof}
	Set $\Omega_{\mathscr{R}}:=B_{\mathscr{R}}\setminus \overline{\Bin}$, where $B_{\mathscr{R}}$ denotes the disc of radius $\mathscr{R}$ centred at zero and $\mathscr{R}>0$ is chosen so large that $\Bin \subset B_{\mathscr{R}}$. We consider the vector field defined in \eqref{eq:pohozaev-vector} for
	$u:=\psi$.
	Note that $u$ is harmonic in $\Bout$, constant on $\cS$, and has the expansion
	\[
		u(x)=\log\abs{x}+O(1)\quad\text{as }\abs{x}\to\infty.
	\]

	Pohozaev's identity \eqref{eq:pohozaev} in dimension two simplifies to
	\[
		\nabla \cdot  X = 0,
	\]
	because $u$ is harmonic in $\Bout$.
	By the divergence theorem on $\Omega_{\mathscr{R}}$ and $\nabla \cdot  X=0$,
	\[
		0=\int_{\Omega_{\mathscr{R}}}\nabla \cdot  X \dd x =\int_{\partial\Omega_{\mathscr{R}}} X\cdot\nu \dd \cH^1
		=\int_{\partial B_{\mathscr{R}}} X\cdot\nu \dd \cH^1 + \int_{\cS} X\cdot\nu \dd \cH^1,
	\]
	where $\nu$ is the outer unit normal of $\Omega_{\mathscr{R}}$.
	On the outer boundary $\partial B_{\mathscr{R}}$, $\nu$ is the radial unit vector $e_r$.
	On the inner boundary $\cS=\partial \Bin$, the outer normal of $\Omega_{\mathscr{R}}$ points into $\Bin$, hence $\nu=-n$, where $n$ is the unit normal of $\cS$ pointing from $\Bin$ into $\Bout$.

	On $\cS$ we have $u=0$, therefore $\nabla u = (\partial_n u)n$.
	Using $\nu=-n$ we compute
	\begin{align*}
		X\cdot\nu & = \big(x\cdot\nabla u\big) \nabla u\cdot\nu - \frac12 \abs{\nabla u}^2  x\cdot\nu \\
		          & = -(x\cdot n)(\partial_n u)^2
		+\frac12(x\cdot n)(\partial_n u)^2                                                            \\
		          & = -\frac12  (x \cdot n) (\partial_n u)^2
	\end{align*}
	on $\cS$.
	Therefore
	\[
		\int_{\cS} X\cdot\nu \dd \cH^1 = -\frac12\int_{\cS} (x \cdot n) (\partial_n u)^2 \dd \cH^1.
	\]

	By \eqref{eq:far-field-derivatives}, on $\partial B_{\mathscr{R}}$ we have
	$x\cdot\nabla u=\mathscr{R}\partial_r u=1+O(1/\mathscr{R})$ and $x\cdot\nu=\mathscr{R}$.
	Hence,
	\[
		\begin{aligned}
			X\cdot\nu
			 & = \big(x\cdot\nabla u\big) \nabla u\cdot\nu - \frac12 \abs{\nabla u}^2  x\cdot\nu                                   \\
			 & = \left(1+O\Big(\frac{1}{\mathscr{R}}\Big)\right)\Big(\frac{1}{\mathscr{R}}+O\Big(\frac{1}{\mathscr{R}^2}\Big)\Big)
			- \frac12 \left(\frac{1}{\mathscr{R}^2}+O\Big(\frac{1}{\mathscr{R}^3}\Big)\right)\mathscr{R}                           \\
			 & = \frac{1}{2\mathscr{R}} + O\Big(\frac{1}{\mathscr{R}^2}\Big).
		\end{aligned}
	\]
	Therefore, we have
	\[
		\int_{\partial B_{\mathscr{R}}} X\cdot\nu \dd \cH^1
		= \int_0^{2\pi}\Big(\frac{1}{2\mathscr{R}}+O\Big(\frac{1}{\mathscr{R}^2}\Big)\Big)\mathscr{R}\dd \theta
		= \pi + O(1/\mathscr{R}) \to \pi \qquad\text{as }\mathscr{R}\to\infty.
	\]
	Collecting the two boundary contributions, we conclude
	\[
		0=\lim_{\mathscr{R}\to\infty}\Big(\int_{\partial B_{\mathscr{R}}} X\cdot\nu \dd \cH^1 + \int_{\cS} X\cdot\nu \dd \cH^1 \Big)
		= \pi - \frac12\int_{\cS} (x\cdot n)  (\partial_n u)^2  \dd \cH^1.
	\]
\end{proof}

We need the following two geometric identities.

\begin{lemma}\label{lem:minkowski}
	We have
	\begin{equation*}
		\int_{\cS} x \cdot n \dd \cH^1 =2\pi,\qquad \int_{\cS} \kappa (x \cdot n) \dd \cH^1 = \cP(\Bin).
	\end{equation*}
\end{lemma}

\begin{proof}
	The first formula is a simple application of the divergence theorem, using that $\abs{\Bin} = \pi$.
	The second is the Minkowski identity for the signed curvature in two dimensions.
	It can be deduced from writing $\cS$ as a curve $\gamma$ in arc-length parametrisation oriented anticlockwise as above and integrating $\frac{\mathrm d}{\mathrm ds}(\gamma(s) \cdot \gamma'(s)) = 1- \kappa(s) (\gamma(s) \cdot n(s))$ over one period.
	Here, the normal vector $n(s)$ is obtained by a clockwise rotation of $\gamma'$ by $\pi/2$, and the signed curvature is defined via the Frenet formula $\gamma''(s) = -\kappa(s) n(s)$.
\end{proof}

For every solution of \eqref{eq:ofbvp}, multiplying the jump condition by
$x\cdot n$, integrating over $\cS$, and applying Lemmas \ref{lem:pohozaev} and
\ref{lem:minkowski} gives
\begin{equation}
	\lambda=\frac{\cP(\Bin)}{2\pi}-\frac{\We}{2}.
	\label{eq:lambda-global-three}
\end{equation}

\begin{proof}[Proof of Proposition \ref{prop:globalWe2}]
	If $\We=0$, the jump condition gives $\kappa=\lambda$. Since $\kappa$ is the signed curvature and $\cS$ is a Jordan curve, the only possibility is a circle. The area constraint gives the unit circle.

	Let us now consider $\We>0$.
	We know that the perimeter $\cP(\Bin) \ge 2 \pi$ by the isoperimetric inequality, since $\abs{\Bin} = \pi$, with equality if and only if $\cS$ is a unit circle.

	Combining \eqref{eq:jumpGB} and \eqref{eq:lambda-global-three}, we obtain
	\begin{equation} \label{eq:fluxL2}
		\int_{\cS} (\partial_n \psi)^2 \dd \cH^1 = \cP(\Bin)-\frac{\cP(\Bin)^2}{\pi\We}+\frac{4\pi}{\We}.
	\end{equation}
	Hence,
	\begin{align*}
		0 & \le \int_{\cS} (\partial_n \psi)^2 \dd \cH^1 - \frac{4\pi^2}{\cP(\Bin)} = \cP(\Bin)-\frac{\cP(\Bin)^2}{\pi\We}+\frac{4\pi}{\We}- \frac{4\pi^2}{\cP(\Bin)} \\
		  & = \frac{1}{\We \pi \cP(\Bin)} (\cP(\Bin)-2\pi)(\cP(\Bin)+2\pi) (\pi \We - \cP(\Bin))
	\end{align*}
	by the Cauchy--Schwarz inequality and \eqref{eq:fluxL1}.
	If $\cS$ is not a unit circle, then $\cP(\Bin)>2\pi$ and we deduce $\We > 2$.
\end{proof}

\begin{remark}
	The estimate uses only the sharp constants in the isoperimetric inequality and the Cauchy--Schwarz inequality.
\end{remark}

\section{Global rigidity of the unit circle for \texorpdfstring{$\We \le 3$}{We <= 3}} \label{sec:globalWe3}
We identify $\R^2$ with $\CC$ and set
\[
	\mathbb D:=\{z\in\CC \; ; \; \abs{z}<1\}.
\]
We first fix the normalisation of the exterior Riemann map used in this section and in Sections \ref{sec:cw} and \ref{sec:var}. Let $E\subset\CC$ be the closure of a bounded Jordan
domain. By the exterior Riemann mapping theorem,
there is a unique conformal map fixing infinity whose leading coefficient at
infinity is positive. We denote this map by $\Phi$ and write its Laurent expansion at infinity as
\begin{equation}\label{eq:exterior-conformal-map}
	\begin{gathered}
		\Phi\colon
		(\CC\cup\{\infty\})\setminus\overline{\mathbb D}
		\longrightarrow(\CC\cup\{\infty\})\setminus E,\\
		\Phi(\infty)=\infty,
		\qquad
		\Phi'(\infty):=
		\lim_{\zeta\to\infty}\frac{\Phi(\zeta)}{\zeta}=c>0,\\
		\Phi(\zeta)
		=c\left(
		\zeta+\beta_0+\sum_{n=1}^{\infty}\beta_n\zeta^{-n}
		\right),
		\qquad \abs{\zeta}>1.
	\end{gathered}
\end{equation}
We call $\Phi$ the normalised exterior conformal map of $E$. The coefficients
$\beta_0,\beta_1,\ldots\in\CC$ are the corresponding normalised Laurent
coefficients.
We next record Gr\"onwall's area formula for the normalised exterior conformal map.

\begin{lemma}\label{lem:exterior-area}
	Let $E\subset\CC$ be compact with boundary given by a rectifiable Jordan curve, and let $\Phi$, $c$, and $(\beta_n)_{n\ge0}$ be as in
	\eqref{eq:exterior-conformal-map}. Then
	$c=\Phi'(\infty)=\capa(E)$ and
	\begin{equation}
		\abs{E}
		=\pi c^2\left(
		1-\sum_{n=1}^{\infty}n\abs{\beta_n}^2
		\right).
		\label{eq:omitted-area}
	\end{equation}
\end{lemma}

\begin{proof}
	The normalisation of the exterior Riemann map gives $c=\capa(E)$; see
	\cite[Thm.~5.2.3]{Ransford1995} and
	\cite[Sec.~2, p.~691]{LiesenSeteNasser2017}. We apply the area computation in
	the proof of Gr\"onwall's area theorem
	\cite[Thm.~1.2.1]{ThomasTuneskiVasudevarao2018} to $G:=c^{-1}\Phi$, for which
	$G'(\infty)=1$, and then rescale by $c$ to obtain
	\eqref{eq:omitted-area}.
\end{proof}

\begin{proof}[Proof of Theorem \ref{thm:globalWe3}]
	Let us assume that $(\cS,\psi,\We,\lambda)$ solves \eqref{eq:ofbvp} and that $\cS$ is not a circle.
	Again, by \cite{KochLeoni} the curve $\cS$ is smooth.

	Following the exterior conformal formulation of \cite[(4)--(7)]{Wegmann2000-ql},
	let $\Phi$ be the normalised exterior conformal map in
	\eqref{eq:exterior-conformal-map} associated with $E=\overline{\Bin}$.
	By Lemma \ref{lem:exterior-area}, its leading coefficient is
	$c=\capa(\overline{\Bin})$.
	Applying the Kellogg--Warschawski theorem to the exterior conformal
	map, after inversion $\zeta=1/z$, we find that $\Phi$ extends smoothly to
	$\abs{\zeta}=1$ and that $\Phi'$ does not vanish there
	\cite[Thms.~3.5 and~3.6, pp.~48--49]{MR1217706}.

	The function $\psi\circ\Phi-\log\abs{\zeta}$ is bounded and harmonic for
	$\abs{\zeta}>1$ and has zero boundary values on $\abs{\zeta}=1$. After the inversion $\zeta=1/z$, the resulting bounded harmonic function
	has a removable singularity at the origin. The maximum principle, followed
	by continuity to the boundary, gives
	\[
		\psi(\Phi(\zeta))=\log\abs{\zeta}
		\qquad\text{for }\abs{\zeta}\ge1.
	\]
	Therefore, on $\zeta=\ee^{\ii t}$, $t\in[0,2\pi)$,
	\begin{equation}\label{eq:conformal-boundary-identities-global-three}
		\dd\cH^1\bigl(\Phi(\ee^{\ii t})\bigr)=\abs{\Phi'(\ee^{\ii t})}\dd t,
		\qquad
		[\partial_n\psi](\Phi(\ee^{\ii t}))=\abs{\Phi'(\ee^{\ii t})}^{-1}.
	\end{equation}
	We set
	\[
		G(z):=\Phi'(1/z)\quad\text{for }0<\abs{z}<1,
		\qquad
		G(0):=c.
	\]
	The function $G$ is analytic and zero-free in the disc and
	extends smoothly and without zeros to the unit circle. Since the disc is
	simply connected, it has an analytic logarithm $L$, which we normalise by
	$L(0)=\log c$; see
	\cite[Ch.~4, Sec.~4.4, Cor.~2, p.~142]{Ahlfors2021}. Let
	$\Lambda=(-\partial_t^2)^{1/2}$ on $2\pi$-periodic
	functions, so that
	\begin{equation} \label{eq:fourier-frac}
		\Lambda \ee^{\ii k t}=\abs{k}\ee^{\ii k t},
		\qquad k \in \Z.
	\end{equation}
	We set $\gamma(t):=\Phi(\ee^{\ii t})$, which parametrises $\cS$
	anticlockwise, the signed curvature is
	\[
		\kappa(\gamma(t))
		=
		\frac{1}{\abs{\gamma'(t)}}
		\operatorname{Im}\left(\frac{\gamma''(t)}{\gamma'(t)}\right).
	\]
	Since
	\[
		\frac{\gamma''(t)}{\gamma'(t)}
		=
		\ii\left(
		1+\ee^{\ii t}\frac{\Phi''(\ee^{\ii t})}{\Phi'(\ee^{\ii t})}
		\right)
	\]
	and $G(\ee^{-\ii t})=\Phi'(\ee^{\ii t})$, analyticity of $L$ and
	\eqref{eq:fourier-frac} give
	\[
		\Lambda\log\abs{\Phi'(\ee^{\ii t})}
		=-\operatorname{Re}\left(\ee^{\ii t}\frac{\Phi''(\ee^{\ii t})}{\Phi'(\ee^{\ii t})}\right).
	\]

	Consequently, the curvature of the conformal parametrisation is
	\begin{equation}\label{eq:curvature-conformal-global-three}
		\kappa(\Phi(\ee^{\ii t}))
		=
		\frac{1-\Lambda\log\abs{\Phi'(\ee^{\ii t})}}{\abs{\Phi'(\ee^{\ii t})}}.
	\end{equation}

	Substituting
	\eqref{eq:conformal-boundary-identities-global-three} and
	\eqref{eq:curvature-conformal-global-three} into the jump condition in
	\eqref{eq:ofbvp}, and using \eqref{eq:lambda-global-three}, we obtain
	\[
		\frac{1-\Lambda\log\abs{\Phi'(\ee^{\ii t})}}{\abs{\Phi'(\ee^{\ii t})}}
		-
		\frac{\We}{2}\abs{\Phi'(\ee^{\ii t})}^{-2}
		=
		\frac{\cP(\Bin)}{2\pi}-\frac{\We}{2}.
	\]
	We set
	\[
		D_\Phi:=\frac{1}{2\pi}\int_0^{2\pi}\abs{\Phi'(\ee^{\ii t})}^2\dd t-1.
	\]
	Multiplying the preceding identity by $\abs{\Phi'(\ee^{\ii t})}^2$ and averaging
	in $t$, using the first equation in
	\eqref{eq:conformal-boundary-identities-global-three}, gives
	\begin{align}
		 & \frac{1}{2\pi}\int_0^{2\pi}
		\abs{\Phi'(\ee^{\ii t})}\,\Lambda\log\abs{\Phi'(\ee^{\ii t})}\dd t
		\label{eq:boundary-energy-global-three} \\
		 & =
		\left(\frac{\We}{2}-\frac{\cP(\Bin)}{2\pi}\right)D_\Phi. \nonumber
	\end{align}

	We now use the analytic square-root representation of the boundary length. Set
	$F:=\exp(L/2)$, so that
	\[
		F(z)^2=G(z),
		\qquad F(0)=c^{1/2},
		\qquad
		F(z)=\sum_{n=0}^{\infty}c_nz^n,
		\qquad z\in\mathbb D.
	\]
	From \eqref{eq:exterior-conformal-map},
	\[
		G(z)
		=c\left(
		1-\sum_{m=1}^{\infty}m\beta_m z^{m+1}
		\right),
	\]
	so the coefficient of $z$ in $F(z)^2$ vanishes. Since $c_0^2=c>0$, we have
	$c_1=0$. The functions $F$ and $L/2=\log F$ have smooth boundary traces.
	Moreover, $G(\ee^{-\ii t})=\Phi'(\ee^{\ii t})$ gives
	\[
		\abs{\Phi'(\ee^{\ii t})}=\abs{F(\ee^{-\ii t})}^2,
		\qquad
		\log\abs{\Phi'(\ee^{\ii t})}
		=2\operatorname{Re}\log F(\ee^{-\ii t}).
	\]
	Changing variables once and using Parseval's identity therefore yields
	\begin{equation}\label{eq:perimeter-F-global-three}
		\frac{\cP(\Bin)}{2\pi}
		=
		\frac{1}{2\pi}\int_0^{2\pi}\abs{\Phi'(\ee^{\ii t})}\dd t
		=
		\frac{1}{2\pi}\int_0^{2\pi}\abs{F(\ee^{\ii t})}^2\dd t
		=
		\sum_{n=0}^{\infty}\abs{c_n}^2.
	\end{equation}
	The same boundary relations give
	\begin{align}
		\frac{1}{2\pi}\int_0^{2\pi}
		\abs{\Phi'(\ee^{\ii t})}\,\Lambda\log\abs{\Phi'(\ee^{\ii t})}\dd t
		 & =
		\frac{1}{2\pi}\int_0^{2\pi}
		\abs{F(\ee^{\ii t})}^2\,2\operatorname{Re}\!\left(\frac{\ee^{\ii t}F'(\ee^{\ii t})}{F(\ee^{\ii t})}\right)\dd t
		\notag \\
		 & =
		2\operatorname{Re}\left(
		\frac{1}{2\pi}\int_0^{2\pi} \ee^{\ii t}F'(\ee^{\ii t})\overline{F(\ee^{\ii t})}\dd t
		\right)
		\notag \\
		 & =
		2\sum_{n=0}^{\infty}n\abs{c_n}^2.
		\label{eq:weighted-F-global-three}
	\end{align}
	By Lemma \ref{lem:exterior-area} and $\abs{\Bin}=\pi$, we have
	\[
		1
		=c^2\left(
		1-\sum_{m=1}^{\infty}m\abs{\beta_m}^2
		\right).
	\]
	On the other hand, Parseval's identity gives
	\[
		\frac{1}{2\pi}\int_0^{2\pi}\abs{\Phi'(\ee^{\ii t})}^2\dd t
		=c^2\left(
		1+\sum_{m=1}^{\infty}m^2\abs{\beta_m}^2
		\right).
	\]
	Since the coefficient of $z^k$ in $F(z)^2=G(z)$ is
	$\sum_{j=0}^{k}c_jc_{k-j}$, comparison with the preceding expansion of $G$
	gives, for $k\ge2$,
	\[
		\sum_{j=0}^{k}c_jc_{k-j}
		=-c(k-1)\beta_{k-1}.
	\]
	Consequently, subtracting the preceding area identity from the Parseval identity
	gives
	\begin{equation}\label{eq:square-gap-global-three}
		D_\Phi=
		\sum_{k=2}^{\infty}\frac{k}{k-1}
		\abs{\sum_{j=0}^{k}c_jc_{k-j}}^2.
	\end{equation}
	In particular, $D_\Phi$ vanishes if and only if $\Phi'$ is constant. By the area normalisation this is equivalent to $\cS$ being a unit circle. Since we assumed that $\cS$ is not a circle we have $D_\Phi>0$.

	We now estimate the right-hand side in \eqref{eq:square-gap-global-three}. As $c_1=0$, for every $k\ge2$ the sum $\sum_{j=0}^{k}c_jc_{k-j}$ has at most $k$ non-zero summands. Therefore Cauchy's inequality, $k/(k-1)\le2$, and \eqref{eq:perimeter-F-global-three} imply
	\begin{align}
		D_\Phi
		 & \le
		2\sum_{k=2}^{\infty}k\sum_{j=0}^{k}\abs{c_j}^2\abs{c_{k-j}}^2
		=
		2\sum_{i,j=0}^{\infty}(i+j)\abs{c_i}^2\abs{c_j}^2
		\notag                      \\
		 & =4\frac{\cP(\Bin)}{2\pi}
		\sum_{j=0}^{\infty}j\abs{c_j}^2.
		\label{eq:hardy-gap-global-three}
	\end{align}

	Now \eqref{eq:boundary-energy-global-three} and \eqref{eq:weighted-F-global-three} imply
	\[
		2\sum_{n=0}^{\infty}n\abs{c_n}^2
		=
		\left(\frac{\We}{2}-\frac{\cP(\Bin)}{2\pi}\right)D_\Phi.
	\]
	We may divide by $D_\Phi$ and obtain
	\[
		\frac{\We}{2}
		=
		\frac{\cP(\Bin)}{2\pi}
		+
		\frac{2\sum_{n=0}^{\infty}n\abs{c_n}^2}{D_\Phi}.
	\]
	Since $D_\Phi>0$, \eqref{eq:hardy-gap-global-three} gives
	\[
		\frac{2\sum_{n=0}^{\infty}n\abs{c_n}^2}{D_\Phi}
		\ge
		\frac{1}{2\frac{\cP(\Bin)}{2\pi}}.
	\]
	Consequently,
	\begin{equation}\label{eq:we-lower-bound-global-three}
		\frac{\We}{2}
		\ge
		\frac{\cP(\Bin)}{2\pi}
		+
		\frac{1}{2\frac{\cP(\Bin)}{2\pi}}.
	\end{equation}
	Since $\abs{\Bin}=\pi$ and $\cS$ is not a unit circle, the isoperimetric inequality gives $\cP(\Bin)>2\pi$. Hence
	\[
		\frac{\cP(\Bin)}{2\pi}
		+
		\frac{1}{2\frac{\cP(\Bin)}{2\pi}}
		-
		\frac{3}{2}
		=
		\frac{\left(\frac{\cP(\Bin)}{2\pi}-1\right)\left(2\frac{\cP(\Bin)}{2\pi}-1\right)}{2\frac{\cP(\Bin)}{2\pi}}
		>
		0.
	\]
	Thus \eqref{eq:we-lower-bound-global-three} implies $\We>3$.
	This proves the theorem.
\end{proof}

\section{The Crowdy--Wegmann branches}
\label{sec:cw}
In this section, we present the non-circular solutions constructed by
Crowdy and Wegmann \cite{Crowdy1999,Wegmann2000-ql} and calculate the parameter
at which their boundary curves first cease to be Jordan curves. We use the same
conformal formulation as in the proof of Theorem
\ref{thm:globalWe3}.

A bubble with boundary prescribed by a Jordan curve determines a normalised exterior conformal map $\Phi$ as in
\eqref{eq:exterior-conformal-map}, and its potential is characterised by
$\psi(\Phi(\zeta))=\log\abs{\zeta}$ for
$\abs{\zeta}>1$. Conversely, suppose that $\Phi$ is holomorphic and univalent
in $\CC\setminus\overline{\mathbb D}$, has the normalisation at infinity from
\eqref{eq:exterior-conformal-map}, and extends $\C^2$-smoothly to the unit
circle, with $\Phi'\ne0$ there and an injective boundary trace. Then
$\Gamma:=\Phi(\partial\mathbb D)$ is a Jordan curve, and the same relation
defines a potential satisfying the first three conditions in
\eqref{eq:ofbvp}. After a dilation which fixes the enclosed area,
\eqref{eq:conformal-boundary-identities-global-three} and
\eqref{eq:curvature-conformal-global-three} express the remaining jump
condition in \eqref{eq:ofbvp} entirely in terms of the boundary values of
$\Phi'$. Thus constructing a bubble reduces to finding such a conformal map
for which this form of the jump condition holds.

Crowdy \cite{Crowdy1999} found the two-fold symmetric solution family,
and Wegmann and Crowdy \cite{Wegmann2000-ql} subsequently generalised it
to $m$-fold symmetric maps for every integer $m\ge2$. Up to translation and
rotation, we may take the parameter $a$ of \cite{Wegmann2000-ql} to be
nonnegative and set $b:=a^m$. The derivative $\Phi_{m,b}'$ is the function
$h_m$ in their equation~(26), and $\Phi_{m,b}$ is the primitive $f_m$ in
their equation~(31). Their exterior map, normalised at infinity as in
\eqref{eq:exterior-conformal-map} but not yet rescaled to enclose area
$\pi$, and its derivative are
\begin{equation}\label{eq:CW-map}
	\Phi_{m,b}(\zeta)
	=
	\zeta+
	\frac{4mb\zeta}{(m-1)(\zeta^m+(m-1)b)},
	\qquad
	\Phi_{m,b}'(\zeta)
	=
	\left(
	\frac{\zeta^m-(m+1)b}{\zeta^m+(m-1)b}
	\right)^2.
\end{equation}
For fixed $m$ and $b$, we set
$\Gamma_{m,b}:=\Phi_{m,b}(\partial\mathbb D)$. We call $b$ admissible if
$\Phi_{m,b}$ is univalent in $\CC\setminus\overline{\mathbb D}$ and its
boundary trace is injective. Then $\Gamma_{m,b}$ is a Jordan curve. As $b$ increases, we are therefore
interested in the first parameter at which two distinct points of
$\partial\mathbb D$ have the same image. At this parameter, two distinct parts
of $\Gamma_{m,b}$ meet, the boundary trace ceases to be injective, and the
curve is no longer a Jordan curve. Thus the first such parameter is the endpoint of the conformal
branch in the class of Jordan curves.

For every admissible parameter, \cite[(35)]{Wegmann2000-ql} gives the area
enclosed by $\Gamma_{m,b}$ as $\pi\mathfrak a_m(b)$, where
\begin{equation}\label{eq:CW-area-factor}
	\mathfrak a_m(b)
	:=
	1-
	\frac{16m^2b^2}{m-1}
	\frac{1+(m-1)b^2}{(1-(m-1)^2b^2)^2}.
\end{equation}
Thus $\mathfrak a_m(b)>0$ for every such parameter. Taking the positive square
root, we set
\[
	\widehat\Phi_{m,b}
	:=\mathfrak a_m(b)^{-1/2}\Phi_{m,b},
	\qquad
	\cS_{m,b}:=\widehat\Phi_{m,b}(\partial\mathbb D).
\]
The curve $\cS_{m,b}$ encloses area $\pi$. On its exterior, we define the
corresponding potential by
\[
	\psi_{m,b}(\widehat\Phi_{m,b}(\zeta))
	=\log\abs{\zeta},
	\qquad \abs{\zeta}>1.
\]
Combining \cite[(29),(30),(35)]{Wegmann2000-ql} with our normalisation
conventions, we find a constant $\lambda_{m,b}$ such that
$(\cS_{m,b},\psi_{m,b},\We_m(b),\lambda_{m,b})$ solves
\eqref{eq:ofbvp}, where
\begin{equation}\label{eq:CW-normalised-Weber-number}
	\We_m(b)
	=
	(m+1)
	\frac{1-(m+1)^2b^2}{1+(m-1)(m+1)b^2}
	\,\mathfrak a_m(b)^{-1/2}.
\end{equation}
The coefficient of $\abs{f'}^{-1}$ in the boundary equation of
\cite{Wegmann2000-ql} is $\We/2$ in our normalisation. At $b=0$, the map
$\Phi_{m,0}$ is the identity and $\We_m(0)=m+1$. Hence the $m$-th branch
emanates from the circular family at the bifurcation value $\We=m+1$. Expanding
\eqref{eq:CW-normalised-Weber-number}, in agreement with
\cite[(39)]{Wegmann2000-ql}, shows that the branches with $m=2,3,4$
initially move towards larger Weber numbers, whereas those with $m\ge5$
initially move towards smaller Weber numbers. In particular,
\[
	\We_2(b)=3+60b^2+O(b^4),
	\qquad b\downarrow0,
\]
and hence the non-circular branch emanating from $\We=3$ immediately lies
in the range $\We>3$.

Wegmann and Crowdy numerically identified these parameters for
$m=2,\ldots,7$ in terms of $a$ and observed that the critical curves undergo an
$m$-fold pinching \cite[p.~2136 and Sec.~5]{Wegmann2000-ql}.
The exact thresholds for $m=2,3,4$ were subsequently obtained by Hayford and
Wang in an equivalent formulation \cite{HayfordWang2022}. We determine the
corresponding thresholds explicitly for every $m\ge2$.
For each integer $k\ge0$, let $U_k$ denote the Chebyshev polynomial of the second kind and
set
\begin{equation}\label{eq:CW-pinching-data}
	D_m(x)
	:=
	-(m-1)U_m(x)+\frac{(m+1)^2}{m-1}U_{m-2}(x),
	\qquad
	M_m:=\max_{-1\le x\le1}\abs{D_m(x)}.
\end{equation}
Using $U_k'(1)=k(k+1)(k+2)/3$, we find that
\[
	D_m(1)=2(m+1),
	\qquad
	D_m(-1)=2(-1)^m(m+1),
	\qquad
	D_m'(1)=-\frac{2m^2(m+1)}{3}<0.
\]
Thus $D_m(x)>D_m(1)$ for $x<1$ sufficiently close to $1$. Consequently,
$M_m>2(m+1)$, and we may define
\begin{equation}\label{eq:CW-pinch-parameter}
	b_m^*
	:=
	\frac{2}{M_m+\sqrt{M_m^2-4(m+1)^2}}.
\end{equation}
In particular, $0<b_m^*<1/(m+1)$.

\begin{lemma}\label{lem:CW-univalence}
	For every $m\ge2$ and $0\le b<b_m^*$, the map $\Phi_{m,b}$ is univalent in
	$\CC\setminus\overline{\mathbb D}$ and its boundary trace is injective.
	At $b=b_m^*$, two distinct points of the unit circle
	have the same image under $\Phi_{m,b}$.
\end{lemma}

\begin{proof}
	Let $z,w\in\partial\mathbb D$ with $z\ne w$. We write
	\[
		z=\ee^{\ii (\varphi+\theta)},
		\qquad
		w=\ee^{\ii (\varphi-\theta)},
		\qquad
		\varphi\in\R,
		\qquad
		0<\theta<\pi.
	\]
	For every integer $k\ge0$, we have
	\[
		\sum_{j=0}^k z^{k-j}w^j
		=
		\ee^{\ii k\varphi}U_k(\cos\theta).
	\]
	Applying this to the difference of the two rational functions
	in \eqref{eq:CW-map}, we obtain the factorisation
	\begin{equation}\label{eq:CW-boundary-factorisation}
		\Phi_{m,b}(z)-\Phi_{m,b}(w)
		=
		\frac{(z-w)\left(
			\ee^{2\ii m\varphi}
			-bD_m(\cos\theta)\ee^{\ii m\varphi}
			+(m+1)^2b^2
			\right)}{
			(z^m+(m-1)b)(w^m+(m-1)b)
		}.
	\end{equation}

	Suppose that $0<b<1/(m+1)$ and that the last factor in the numerator
	vanishes. Setting $\omega:=\ee^{\ii m\varphi}$, dividing \eqref{eq:CW-boundary-factorisation} by $\omega$, and
	taking imaginary parts gives
	\[
		\bigl(1-(m+1)^2b^2\bigr)\sin(m\varphi)=0.
	\]
	Thus $\omega\in\{-1,1\}$, and the real part of \eqref{eq:CW-boundary-factorisation} yields
	\begin{equation}\label{eq:CW-pinching-criterion}
		\abs{D_m(\cos\theta)}
		=
		\frac1b+(m+1)^2b.
	\end{equation}
	Conversely, if \eqref{eq:CW-pinching-criterion} holds for some
	$0<\theta<\pi$, we can choose $\varphi$ such that
	$\ee^{\ii m\varphi}$ has the sign of $D_m(\cos\theta)$. Then
	\eqref{eq:CW-boundary-factorisation} shows that two distinct
	boundary points have the same image under $\Phi_{m,b}$.

	The function $b\mapsto b^{-1}+(m+1)^2b$ is strictly decreasing on
	$(0,1/(m+1))$. Moreover,
	\eqref{eq:CW-pinch-parameter} is its unique solution of
	\[
		\frac1b+(m+1)^2b=M_m
	\]
	in this interval. Hence \eqref{eq:CW-pinching-criterion} cannot hold
	for $0<b<b_m^*$. Since $M_m>2(m+1)$, every maximiser of $\abs{D_m}$ lies in
	$(-1,1)$. Choosing such a maximiser and applying the converse implication
	above, we obtain distinct points $z,w\in\partial\mathbb D$ such that
	\[
		\Phi_{m,b_m^*}(z)=\Phi_{m,b_m^*}(w).
	\]
	Finally, the formulae in \eqref{eq:CW-map} show that $\Phi_{m,b}$ has
	neither a pole nor a critical point in
	$\CC\setminus\mathbb D$ for $0\le b<b_m^*$.
	Its boundary trace is injective
	by the preceding argument, and
	$\Phi_{m,b}(\zeta)=\zeta+O(\abs{\zeta}^{1-m})$ as $\abs{\zeta}\to\infty$. The argument principle
	\cite[Ch.~4, Sec.~5.2]{Ahlfors2021} now shows that
	$\Phi_{m,b}$ maps $\CC\setminus\overline{\mathbb D}$ conformally onto the unbounded component
	of the complement of its boundary curve. This proves the claim.
\end{proof}

\begin{proposition}\label{prop:CW-large-Weber}
	For every $\We\ge10$, there exists a non-circular Crowdy--Wegmann
	bubble of area $\pi$ which solves \eqref{eq:ofbvp} with Weber number
	$\We$.
\end{proposition}

\begin{proof}
	For $0<\theta<\pi$, the trigonometric representation of the Chebyshev
	polynomials gives
	\begin{equation}\label{eq:CW-D-trigonometric}
		D_m(\cos\theta)
		=
		\frac{2}{m-1}\left(
		2m\cot\theta\sin(m\theta)
		-(m^2+1)\cos(m\theta)
		\right).
	\end{equation}
	The values at $\theta=0$ and $\theta=\pi$ are understood by continuous
	extension.
	We claim that
	\[
		M_m\le3(m+1)
		\qquad m\ge4.
	\]
	At the endpoints, the identities above give
	\[
		\abs{D_m(1)}
		=
		\abs{D_m(-1)}
		=
		2(m+1)
		\le 3(m+1).
	\]
	Since $M_m>2(m+1)$, the maximum is attained at some
	$\theta\in(0,\pi)$ with
	\[
		\abs{D_m(\cos\theta)}=M_m.
	\]
	Differentiating \eqref{eq:CW-D-trigonometric} at this maximiser shows that
	the desired estimate is equivalent, with $y=\sin^2\theta$ and
	$h=m^2-1$, to
	\[
		5h^2y^2-(16h+36)y+20\ge0.
	\]
	The discriminant of this quadratic polynomial is
	$-144(h-9)(h+1)<0$ for $m\ge4$, while its leading coefficient is
	positive. This proves the claim. Recalling
	\eqref{eq:CW-pinch-parameter}, we deduce that
	\[
		b_m^*
		\ge
		\frac{3-\sqrt5}{2(m+1)}
		>\frac{1}{3(m+1)}.
	\]

	Let $m \in \N $ with $m\ge10$ and choose $b_0:=1/(3(m+1))$. Lemma
	\ref{lem:CW-univalence} shows that $\Phi_{m,b}$ is univalent for
	$0\le b\le b_0$. Substituting $b_0$ into
	\eqref{eq:CW-normalised-Weber-number}, we obtain
	\begin{equation}\label{eq:CW-large-Weber-test-point}
		\bigl[\We_m(b_0)\bigr]^2
		=
		\frac{
			16(m-1)(m+1)^2(2m+1)^2(m+2)^2
		}{
			(5m+4)^2(4m^3-m^2-8m-4)
		}
		<m^2.
	\end{equation}
	The last inequality is equivalent to $Q(m)>0$, where we set
	\[
		Q(m)
		:=
		36m^7-249m^6-960m^5-900m^4
		+176m^3+720m^2+384m+64.
	\]
	Writing $m=n+10$ with an integer $n\ge0$, every coefficient of $Q(n+10)$ is
	strictly positive.
	Hence $Q(m)>0$ for $m\ge10$. Moreover,
	$0<b_0<1/(m+1)$, and therefore
	\eqref{eq:CW-normalised-Weber-number} gives $\We_m(b_0)>0$. In view
	of \eqref{eq:CW-large-Weber-test-point}, we conclude that
	$\We_m(b_0)<m$.

	As noted above, $\We_m(0)=m+1$. Hence continuity yields every value in
	$[m,m+1)$ for some $0<b<b_m^*$. Since
	$\bigcup_{m\ge10}[m,m+1)=[10,\infty)$, the proposition follows.
\end{proof}

\begin{remark}\label{rem:noncircular}
	A Mathematica investigation of the branch ranges in $3<\We<10$ suggests that the
	non-circular Crowdy--Wegmann branches cover all Weber numbers $\We>3$,
	except possibly in two intervals. Their endpoint pairs are
	estimated to be near
	\[
		8.18307625\ldots,\ 8.37519009\ldots
		\qquad\text{and}\qquad
		9,\ 9.06929546\ldots.
	\]
\end{remark}

\section{On a related variational problem} \label{sec:var}
We now turn to the variational problem \eqref{eq:infconvF}.

\subsection{Minimality of the disc for \texorpdfstring{$0 \le \We \le 3$}{0 <= We <= 3}}

For $1\le p<\infty$ and an analytic function $u$ on $\mathbb D$, we use the
Hardy norm
\[
	\norm{u}_{\Hardy^p}^p
	:=
	\sup_{0<r<1}\frac1{2\pi}\int_0^{2\pi}\abs{u(r\ee^{\ii t})}^p\dd t.
\]
We also write
\[
	\norm{u}_{A^2}^2:=\frac1\pi\int_{\mathbb D}\abs{u(z)}^2\dd A(z)
\]
for the Bergman norm, where $\dd A$ denotes planar Lebesgue
measure. We denote by $\Hardy^p(\mathbb D)$ and $A^2(\mathbb D)$ the corresponding
spaces of analytic functions with finite norm.
For $u\in\Hardy^p(\mathbb D)$, we write
\[
	u^*(\ee^{\ii t})
	:=
	\lim_{r\uparrow1}u(r\ee^{\ii t})
\]
for its radial boundary value, which exists for Lebesgue-almost every
$t\in[0,2\pi)$.
For further background on Hardy and Bergman spaces, see
\cite[Ch.~5 and Sec.~8.1]{MR2109650}.

The key tool for Theorem \ref{thm:exact-threshold} is the following
area-constrained perimeter--capacity estimate, whose proof uses three auxiliary
lemmas.

\begin{proposition}
	\label{prop:critical-three}
	Let $E\subset\R^2\simeq\CC$ be compact with boundary given by a rectifiable Jordan curve, and suppose that $\abs{E}=\pi$. Then
	\[
		\cP(E)-2\pi
		\ge
		2\pi\bigl(\sqrt{5\capa(E)^2-1}-\capa(E)-1\bigr).
	\]
	Moreover, equality holds if and only if $E$ is a translate of the unit disc.
\end{proposition}

The first one is the classical area--capacity inequality; see
\cite[Thm.~5.3.5]{Ransford1995}.

\begin{lemma}\label{lem:unit-area-capacity}
	Let $E$ satisfy the assumptions of Proposition \ref{prop:critical-three}.
	Then
	\[
		\capa(E)\ge1.
	\]
\end{lemma}

\begin{proof}
	Let $\Phi$, $c$, and $(\beta_n)_{n\ge0}$ be as in
	\eqref{eq:exterior-conformal-map}. By Lemma \ref{lem:exterior-area},
	$c=\capa(E)$.
	Since $\abs{E}=\pi$, equation \eqref{eq:omitted-area} gives
	\[
		1
		=c^2\left(
		1-\sum_{n=1}^{\infty}n\abs{\beta_n}^2
		\right)
		\le c^2. \qedhere
	\]
\end{proof}

Next, we express the perimeter as a Hardy norm.

\begin{lemma}\label{lem:perimeter-hardy}
	Let $E$ satisfy the assumptions of Proposition
	\ref{prop:critical-three}, and suppose that
	$0\in\operatorname{int}(E)$. Let $\Phi$ be the corresponding normalised
	exterior conformal map, with $c$ and $(\beta_n)_{n\ge0}$ as in
	\eqref{eq:exterior-conformal-map}.
	Define
	\[
		q(z):=\sum_{n=1}^{\infty}n\beta_nz^{n+1},
		\qquad z\in\mathbb D.
	\]
	Then
	\[
		q(0)=q'(0)=0,
		\qquad
		1-q(z)\neq0\quad\text{for all }z\in\mathbb D.
	\]
	Here and below, $q/z^2$ denotes its analytic continuation at the origin,
	with value $q''(0)/2$ there.
	The area and perimeter satisfy
	\begin{equation}
		\norm{\frac{q}{z^2}}_{A^2}^2
		=
		\sum_{n=1}^{\infty}n\abs{\beta_n}^2
		=
		1-\frac{\abs{E}}{\pi c^2}.
		\label{eq:area-conformal-map}
	\end{equation}
	\begin{equation}\label{eq:perimeter-conformal-map}
		\frac{\cP(E)}{2\pi c}=\norm{1-q}_{\Hardy^1}.
	\end{equation}
\end{lemma}

\begin{proof}
	The coefficient formula for the Bergman norm
	\cite[proof of Thm.~5.3.4, p.~77]{MR2109650}, together with
	\eqref{eq:omitted-area}, gives
	\eqref{eq:area-conformal-map}.
	By Carath\'eodory's theorem, $\Phi$ extends continuously and one-to-one to
	$\partial\mathbb D$; see \cite[Ch.~6, Sec.~1.2]{Ahlfors2021}.
	For $0<\abs{z}<1$, set
	\[
		f(z):=\frac1{\Phi(1/z)},
		\qquad
		h(z):=z\Phi(1/z),
	\]
	and set $f(0):=0$ and $h(0):=c$. The Laurent expansion of $\Phi$ shows
	that these are the removable analytic extensions to $\mathbb D$.
	Inversion $z \mapsto \frac1z$ is bi-Lipschitz in a neighbourhood of
	$\partial E$. Hence $f$ maps $\mathbb D$ conformally onto a bounded Jordan
	domain with rectifiable boundary. By the rectifiable-boundary theorem
	\cite[Thm.~3.12, p.~44]{Duren1970}, we have $f'\in\Hardy^1(\mathbb D)$. The same
	boundary extension shows that $h$ is bounded. Therefore, for
	$0<\abs{z}<1$,
	\[
		1-q(z)=\frac1c\Phi'(1/z)=\frac1c h(z)^2f'(z).
	\]
	The identity extends analytically across $0$; since $h$ is bounded and
	$f'\in\Hardy^1(\mathbb D)$, it follows that
	$1-q\in\Hardy^1(\mathbb D)$.
	The boundary map $t\mapsto f(\ee^{\ii t})$ is absolutely continuous and satisfies
	\[
		\frac{\mathrm d}{\mathrm dt}f(\ee^{\ii t})=\ii\ee^{\ii t}(f')^*(\ee^{\ii t})
	\]
	for almost every $t$; see \cite[Thm.~3.11, p.~42]{Duren1970}. Since $f$ is
	bounded away from zero on $\partial\mathbb D$, the boundary parametrisation
	$\gamma(t):=\Phi(\ee^{\ii t})=1/f(\ee^{-\ii t})$ is absolutely continuous and satisfies
	\[
		\abs{\gamma'(t)}
		=
		c\abs{(1-q)^*(\ee^{-\ii t})}
		\qquad\text{for almost every }t\in[0,2\pi].
	\]
	The boundary-value identity for the Hardy norm
	\cite[Thm.~5.1.8]{MR2109650} now gives
	$\cP(E)=2\pi c\norm{1-q}_{\Hardy^1}$.
	Finally, $1-q(0)=1$, while for $0<\abs{z}<1$ it is nonzero because $\Phi'(1/z)\neq0$.
	This proves the claim.
\end{proof}

It remains to relate the Hardy norm in
\eqref{eq:perimeter-conformal-map} to the Bergman norm in
\eqref{eq:area-conformal-map}. The estimate needed here
is the following tailored consequence of Burbea's sharp product inequality
\cite[Thm.~4.1, p.~257]{MR882113}.

\begin{lemma}\label{lem:hardy-bergman-estimate}
	Let $q$ be analytic in $\mathbb D$ and suppose that
	\[
		q(0)=q'(0)=0,
		\qquad
		1-q\in\Hardy^1(\mathbb D),
		\qquad
		1-q(z)\neq0\quad\text{for all }z\in\mathbb D.
	\]
	Then $q/z^2\in A^2(\mathbb D)$,
	\[
		M:=\norm{1-q}_{\Hardy^1}\ge1,
	\]
	and
	\begin{equation}\label{eq:hardy-bergman-estimate}
		\norm{\frac{q}{z^2}}_{A^2}^2\le(M-1)(M+3).
	\end{equation}
	Moreover, equality holds if and only if $q\equiv0$.
\end{lemma}

\begin{proof}
	Since $\mathbb D$ is simply connected and $1-q$ is zero-free, it has an
	analytic square root $g$, normalised by $g(0)=1$. Thus
	\[
		g^2=1-q.
	\]
	Then $g\in\Hardy^2(\mathbb D)$ and
	\[
		\norm{g}_{\Hardy^2}^2=\norm{1-q}_{\Hardy^1}=M.
	\]
	In particular, $M\ge\abs{g(0)}^2=1$. Since $q'(0)=0$, we have
	$g'(0)=0$, and hence
	\[
		g=1+z^2v
	\]
	for some $v\in\Hardy^2(\mathbb D)$. Parseval's identity
	\cite[Sec.~3.4, p.~43]{MR2109650} gives
	\[
		\norm{v}_{\Hardy^2}^2=M-1,
		\qquad
		\norm{g+1}_{\Hardy^2}^2=M+3.
	\]
	Moreover,
	\[
		\frac{q}{z^2}=-(g+1)v.
	\]

	Applying Burbea's product inequality
	\cite[Thm.~4.1, p.~257]{MR882113} gives
	\[
		\norm{\frac{q}{z^2}}_{A^2}^2
		\le
		\norm{g+1}_{\Hardy^2}^2\norm{v}_{\Hardy^2}^2
		=
		(M+3)(M-1).
	\]
	If $v\not\equiv0$, the equality statement in Burbea's product inequality
	\cite[Thm.~4.1, p.~257]{MR882113} gives constants
	$C_1,C_2\in\CC\setminus\{0\}$ and $\zeta\in\mathbb D$ such that
	\[
		g(z)+1=\frac{C_1}{1-\overline\zeta z},
		\qquad
		v(z)=\frac{C_2}{1-\overline\zeta z}.
	\]
	The linear coefficient of $g+1$ vanishes, while its constant coefficient
	is equal to $2$. Thus $\zeta=0$, so both functions are constant. The identity
	$g+1=2+z^2v$ then gives $v\equiv0$, a contradiction. Consequently,
	equality holds only if $v\equiv0$, which is equivalent to $q\equiv0$.
	The converse is immediate.
\end{proof}

\begin{proof}[Proof of Proposition \ref{prop:critical-three}]
	By translation invariance, we may suppose that $0\in\operatorname{int}(E)$.
	We use the notation of \eqref{eq:exterior-conformal-map} and Lemma
	\ref{lem:perimeter-hardy} and set
	\[
		c:=\capa(E),
		\qquad
		\eta:=\norm{\frac{q}{z^2}}_{A^2}^2,
		\qquad
		M:=\norm{1-q}_{\Hardy^1}.
	\]
	Since $\abs{E}=\pi$, equations \eqref{eq:area-conformal-map} and
	\eqref{eq:perimeter-conformal-map} give
	\[
		\eta=1-c^{-2},
		\qquad
		M=\frac{\cP(E)}{2\pi c}.
	\]
	Applying Lemma \ref{lem:hardy-bergman-estimate},
	we obtain
	\[
		\eta\le(M-1)(M+3).
	\]
	Since $M\ge1$, this is equivalent to
	\[
		M\ge-1+\sqrt{4+\eta}.
	\]
	Consequently,
	\[
		\frac{\cP(E)}{2\pi}
		=cM
		\ge
		c\left(-1+\sqrt{5-c^{-2}}\right)
		=
		\sqrt{5c^2-1}-c,
	\]
	which proves the estimate.

	If equality holds in the estimate of Proposition
	\ref{prop:critical-three}, then equality holds in
	\eqref{eq:hardy-bergman-estimate}. Lemma
	\ref{lem:hardy-bergman-estimate} gives $q\equiv0$. Hence $\beta_n=0$ for every
	$n\ge1$, while \eqref{eq:area-conformal-map} gives $c=1$. Thus
	$\Phi(\zeta)=\zeta+\beta_0$, and $E$ is a translate of the unit disc. Conversely,
	a translate of the unit disc gives equality.
\end{proof}

The following logarithmic form is adapted to the variational functional.

\begin{lemma}\label{lem:critical-three-log}
	Let $E$ satisfy the assumptions of Proposition \ref{prop:critical-three}.
	Then
	\[
		\cP(E)-2\pi\ge3\pi\log\capa(E).
	\]
	Moreover, equality holds if and only if $E$ is a translate of the unit disc.
\end{lemma}

\begin{proof}
	Set $c:=\capa(E)$. Lemma \ref{lem:unit-area-capacity} gives $c\ge1$.
	In view of Proposition
	\ref{prop:critical-three}, it remains to show that
	\[
		\sqrt{5c^2-1}-c\ge1+\frac32\log c
		\qquad\text{for }c\ge1.
	\]
	Define $F\colon[1,\infty)\to\R$ by
	\[
		F(c):=\sqrt{5c^2-1}-c-1-\frac32\log c.
	\]
	Then $F(1)=0$ and
	\[
		F'(c)=\frac{5c}{\sqrt{5c^2-1}}-1-\frac{3}{2c}.
	\]
	Rewriting the numerator gives
	\[
		F'(c)
		=
		\frac{(c-1)(80c^3+20c^2-21c-9)}
		{2c\sqrt{5c^2-1}\left(10c^2+(2c+3)\sqrt{5c^2-1}\right)}.
	\]
	The cubic factor is positive on $[1,\infty)$ because
	\[
		80c^3+20c^2-21c-9
		=70+(c-1)(80c^2+100c+79).
	\]
	Thus $F'(c)\ge0$ for
	$c\ge1$, with strict inequality for $c>1$. This proves the estimate.

	If $\cP(E)-2\pi=3\pi\log\capa(E)$, then $c=1$ and $\cP(E)=2\pi$. Hence equality also holds
	in Proposition \ref{prop:critical-three}, and $E$ is a translate of the unit
	disc. The converse is immediate.
\end{proof}

\begin{proposition}\label{prop:disk-minimiser-up-to-three}
	For every $0\le\We\le3$, the unique minimisers of $\F_{\We}$ among compact sets whose boundary is given by a rectifiable Jordan curve and area $\pi$ are unit discs up to translation.
\end{proposition}

\begin{proof}
	Let $E$ be admissible and set $c:=\capa(E)$.
	By Lemma \ref{lem:critical-three-log},
	\[
		\cP(E)-3\pi\log c\ge2\pi.
	\]
	Lemma \ref{lem:unit-area-capacity} gives $c\ge1$.
	Hence, for $0\le\We\le3$,
	\[
		\begin{aligned}
			\F_{\We}(E)
			 & =\cP(E)-\We\pi\log c                            \\
			 & =\bigl(\cP(E)-3\pi\log c\bigr)+(3-\We)\pi\log c \\
			 & \ge2\pi.
		\end{aligned}
	\]
	For the closed unit disc, $\cP(\overline{\mathbb D})=2\pi$ and
	$\capa(\overline{\mathbb D})=1$; see
	\cite[Table~5.1]{Ransford1995}. Thus $\F_{\We}(\overline{\mathbb D})=2\pi$.

	If $\F_{\We}(E)=2\pi$ for an admissible set $E$, then
	$\cP(E)-3\pi\log c=2\pi$. Lemma \ref{lem:critical-three-log} implies that
	$E=x_0+\overline{\mathbb D}$ for some $x_0\in\R^2$. This completes the proof.
\end{proof}

\subsection{For \texorpdfstring{$\We>3$}{We>3}, elliptic deformations lower the energy}

Let us now prove that the threshold is sharp. We consider the one-parameter
family of closed ellipses
\[
	\E_t:=\{(x_1,x_2)\in\R^2 \; ; \;
	\ee^{-2t}x_1^2+\ee^{2t}x_2^2\le1\},
	\qquad t\in\R.
\]
Then $\E_0=\overline{\mathbb D}$, and each $\E_t$ is a convex and compact set
satisfying $\abs{\E_t}=\pi$.

\begin{proposition}\label{prop:localvar}
	For every $\We>3$, there exists $\delta_{\We}>0$ such that
	\[
		\F_{\We}(\E_t)<\F_{\We}(\overline{\mathbb D})
		\qquad\text{whenever}\qquad
		0<\abs{t}<\delta_{\We}.
	\]
\end{proposition}

\begin{proof}
	The logarithmic capacity of an ellipse equals half the sum of its
	semiaxes; see \cite[Table~5.1]{Ransford1995}. We therefore obtain
	\begin{align*}
		\cI(\E_t)
		 & =-\log\capa(\E_t)
		=-\log\!\left(\frac{\ee^t+\ee^{-t}}{2}\right)
		=-\log(\cosh t)         \\
		 & =-\frac12t^2+O(t^4),
		\qquad t\to0.
	\end{align*}
	The perimeter of $\E_t$ can be written as
	\begin{equation}\label{eq:P-integral}
		\cP(\E_t)
		=\int_0^{2\pi}
		\sqrt{\ee^{2t}\sin^2\theta+\ee^{-2t}\cos^2\theta}
		\dd\theta,
	\end{equation}
	from which we obtain the expansion
	\[
		\cP(\E_t)
		=2\pi+\frac{3\pi}{2}t^2+O(t^4),
		\qquad t\to0.
	\]
	Consequently,
	\[
		\F_{\We}(\E_t)-\F_{\We}(\overline{\mathbb D})
		=-\frac{\pi}{2}(\We-3)t^2+O_{\We}(t^4).
	\]
	For $\We>3$, this is negative for all sufficiently small nonzero $t$.
\end{proof}

\subsection{Conclusion of the sharp threshold}

\begin{proof}[Proof of Theorem \ref{thm:exact-threshold}]
	Every planar convex and compact set of positive area has boundary prescribed by a rectifiable Jordan curve. Hence Proposition
	\ref{prop:disk-minimiser-up-to-three} applies to the admissible class in
	\eqref{eq:infconvF} and proves the claim for $0\le\We\le3$. For $\We>3$,
	Proposition \ref{prop:localvar} supplies an admissible ellipse with lower
	energy than the disc.
\end{proof}

\begin{remark}
	For $\We>3$, it remains open to characterise minimisers of
	\eqref{eq:infconvF} and, more broadly, those among closures of bounded Jordan
	domains of area $\pi$.
	It is also an interesting open problem whether, besides the disc and the
	bifurcation branches constructed in \cite{Wegmann2000-ql}, $\F_{\We}$ has
	critical sets whose boundaries, together with some $\psi$ and $\lambda$,
	solve \eqref{eq:ofbvp}. For $\We \gg 3$, numerical experiments yield
	non-circular elongated shapes reminiscent of a stadium.
\end{remark}

\vfill

\end{document}